\def\a{\alpha}
\def\b{\beta}
\def\g{\gamma}
\def\d{\delta}

\def\ve{\varepsilon}

\def\s{\sigma}

\def\vp{\varphi}

\def\es{\emptyset}

\def\sm{\setminus}

\def\sb{\subset}      \def\sbe{\subseteq}

\def\tR{\mathchoice{\mbox{\bf R}}{\mbox{\bf R}}{\mbox{\scriptsize \bf
R}}{\mbox{\tiny\bf R}}}

\def\tN{\mathchoice{\mbox{\bf N}}{\mbox{\bf N}}{\mbox{\scriptsize \bf
N}}{\mbox{\tiny\bf N}}}

\def\tZ{\mathchoice{\mbox{\bf Z}}{\mbox{\bf Z}}{\mbox{\scriptsize \bf
Z}}{\mbox{\tiny\bf Z}}}

\def\bydef{\,\lower-.1ex\hbox{\rm :}\!=}

\def\cA{{\cal A}}
\def\cB{{\cal B}}
\def\cC{{\cal C}}
\def\cD{{\cal D}}
\def\cE{{\cal E}}
\def\cF{{\cal F}}

\def\cP{{\cal P}}

\def\dist{\mbox{ \rm dist} }
\def\Id{\mbox{ \rm Id} }

\def\ol{\overline}

\documentstyle[12pt,titlepage]{article}
\def\bs{\vskip 0.3cm}
\def\proof{\goodbreak\noindent{\sc Proof. }\nobreak}

\def\rem{\vskip3pt\noindent{\sc Remark.} }
\def\endproof{\par\nobreak\hbox to \hsize{\hfil\vrule width 5pt height
5pt}\goodbreak\vskip 3pt}
\newtheorem{theor}{Theorem}

\newtheorem{lemma}[theor]{Lemma}

\def\supp{\mbox{\rm supp}}

\def\Lip{\mbox{\rm Lip}}
\def\diam{\mbox{ \rm diam }}
\def\supp{\mbox{ \rm supp }}

\def\span{\mbox{ \rm  span }}
\def\1{\mbox{\bf 1}}

\title{The Banach space $H^1(X,d,\mu)$, II}
\author{Paul F.X. M\"uller\\ Department of Mathematics\\ J. Kepler
University\\ A-4040 Linz, Austria
Keywords, Hardy space, Isomorphic classification}
\date{}

\titlepage
\begin{document}
\parindent0pt
\maketitle
\newpage
\pagestyle{plain}
\pagenumbering{arabic}

\section{Introduction}

In this paper we give the isomorphic classification of atomic $H^1(X,d,\mu)$,
where $(X,d,\mu)$ is a space of homogeneous type, hereby completing a line of
investigation opend by the word of Bernard Maurey [Ma1], [Ma2], [Ma3] and
continued by Lennard Calrleson [C] and Przemyslaw Wojtaszczyk [Woj1], [Wpj2].

The resulting isomorphic representatives dyadic $H^1, (\sum H^1_n)_{l^1}$ and
$l^1$; each isomorphic type is characterized by geometric properties of
$(X,d,\mu)$.\vskip 0.3cm
Technically, the present paper deals with the existence of \lq\lq Franklin-type"
function on a space of homogeneous type $(X,d,\mu)$ and with the boundedness (in
$H^1(X,d,\mu)$) of the orthogonal projection onto the span of this \lq\lq
Franklin-system".

To this end a construction of S. Jaffard and Y. Meyer [J-M], together with
methods of P. Wojtaszczyk [Woj2] have been adapted. (See Sections 3,4).

Applying Wojtaszczyk's method one concludes that the \lq\lq Franklin-system" in
$H^1(X,d,\mu)$ is equivalent to special three valued martingale differences in
some martingale $H^1$ space; (see Section 4).

Using repeateally Pe\l czy\'nski's decomposition method, this allows us to reduce
in Section 5 the classification problem for atomic $H^1$ spaces to some of the
other's precious work on this subject:
\begin{itemize}
\item Classification of (the span of) special
three valued martingale differences in martingale $H^1$ space [M\"u2].
\item Classification of the isomorphic type of martingale $H^1$ space [M\"u1].
\item Atomic $H^1$ spaces are isomorphic to complemented subspaces of martingale
$H^1$ [M\"u3].\end{itemize}
\section{Preliminaries concerning spaces of homogeneous type}

{\bf Definition and Notation}

Let $d:X\times X\to\tR^+$ be a quasi metric on a set $X$ and let $B(x,r):=\{y\in
X:d(x,y)<r\}$.
Let $\mu$ be a non negative measure on $X$.

Let $A_1, A_2$, $K_1,K_2$ $(K_2\le 1\le K_1)$ be positive finite constants such
that for each $x\in X$ and $r>0$ the following relations hold:

\begin{tabular}{rll}
(1)&$A_1r\le \mu(B(x,r))$\quad\quad\quad &if $r\le K_1\mu(X)$\\
(2)&$B(x,r)=X$&if $r>K_1\mu(X)$\\
(3)&$A_2r\ge\mu(B(x,r))$&if $r\ge K_2\mu(\{x\})$\\
(4)&$B(x,r)=\{x\}$& if $r<K_2\mu(\{x\})$.\end{tabular}

We assume mereover that there exists $\a>0$ and $K_0>0$ such that for each
$x,y,x\in X$
$$|d(x,y)-d(y,z)|\le K_0r^{1-\a}d(x,y)^\a\leqno(5)$$
whenever $d(x,z)<r$ and $d(y,z)<r$.

Following Masias \& Segovia [M-S1,2] a set $X$ equipped with a
measure $\mu$ and a quasi metric $d$ satisfying
$(1)-(5)$ is called a \lq\lq normal space of order $\a$".   The standard
reference to there spaces is the article by Coifman \& Weiss [C-W].

On $(X,d,\mu)$ certain analogues of dyadic intervals have been constructed by G.
David and M. Christ; see  [Dv] and [Ch].

\begin{theor} There exists a collection $\{Q^k_i\sb X:k\in\tZ,i\in I_k\}$ and
constants $\d\in(0,1)$, $\a_0>0, \eta>0$ and $C_1,C_2<\infty$ such that:

 $\mu(X\sm\bigcup Q_i^k)=0$.

If $l\ge k$ then either $Q^l_j\sb Q^k_i$ or $Q_{\rho_j}^l\cap Q^k_i=0$.

For each $(k,i)$ and $l<k$ there exists a unique $j$ such that
$Q^k_i\sb Q^l_j$.

{\rm diam} $Q^k_i\le C_1\d^k$.

Each $Q^k_i$ contains a ball $B(z^k_i,a_0,\d^k)$.

$\mu\{x\in Q^k_i:(x,X\sm Q^k_i)\le td^k\}\le C_2t^\eta\mu(Q^k_i)$
for each $k\in\tZ$, $i\in I_k$ and $t\ge 0$.
                      \end{theor}

As we shall see, the structure 
of this collection determines 
the isomorphic type
of $H^1(X,d,\mu)$. 
However we have 
to discard measure and diameter is too big.

\begin{lemma} In every normal 
space of order $\a$ $(X,d,\mu)$ there exists $L>0$
depending on so that for every nonempty $Q\sb X$ we
have: $\mu(Q)/\diam Q>L$ implies $Q$ consists of exactly one point.\end{lemma}

\proof Select $M>1$ so that $_2M>1$. Then consider two cases.

{\bf Case 1} Suppose there exists $x_0\in Q$ so that $\mu(\{x_0\})\ge
M\diam Q$ then $\{x:d(x,x_0)<KM\diam Q\}=\{x_0\}$. Hence $Q=\{x_0\}$.

{\bf Case 2} For each $x_0\in Q$  we have $\mu(\{x_0\})\le M\diam Q$;
then as $MK_2>1$

\begin{eqnarray*}
\mu(Q)& \le & \mu(B(x_0, K_2M\diam Q)\\
&\le& A_2MK_2\diam Q.\end{eqnarray*}

Let now $L=A_2MK_2$ then we have either $\mu(Q)/\diam Q\le L$ or $Q=\{x_0\}$ for
some $x_0\in X$.
\endproof

Let $\cE:=\{Q^k_\a:k\in\tZ,\a\in I_k\}$. And let $\cF_n$ be the $\s$-Algebra
generated the $n$-th generation of $\cE$ and $\cF_{n-1}$.

The following properties of $\cF_n$ are easily observed:

\begin{enumerate}
\item There exits $N_0\in\tN$, depending on $\d$ (and the geometry of
$(X,d,\mu)$) so that for every $Q\in\cE$ the cardinality of $G_1(Q|\cE)$ is
bounded by $N_0$.
\item There exists $L_0$, depending on $\d$, (and the geometry of $(X,d,\mu)$)
so that for every $Q\in\cE$ and every $P\in G_1(Q|\cE)$ we have
$$\frac{\mu(P)}{\mu(Q)}\ge \frac 1{L_0}.$$
\item Moreover for $Q\in\cE$ we have $\frac 1C\diam Q\le \mu(Q)\le C\diam Q$
where $C$ is as in Lemma 1.\end{enumerate}

The collection $\cE$ has been linked to problems conscerning the isomorphic
structure of $H^1(X,d,\mu)$; see [M\"u].

There we found finitely many sequences of incereasing, pure by atomic
$\a$-algebras
$$[\cF_{n^1}]^\infty_{n=1},\dots,[\cF^N_n]^\infty_{n=1}$$
so that $H^1(X,f,\mu)$ is isomorphic to a complemented subspace of the direct
sum of the related martingale $H^1$-spaces, namely to
$$H^1([\cF^1_n]^\infty_{n=1})\oplus\dots\oplus H^1([\cF^N_n]^\infty _{n=1}).$$
Although not stated explicitely there, when combined with the result in [M\"u]
one observes the following implications:
\begin{enumerate}
\item If $$E=\{t\in X:t\mbox{ lies in infinitely many } Q\in\cC\}$$
satisfies $\mu(E)=0$, then for each $j\le N$:
$H^1([\cF^j_n]^\infty_{n=1})$ is isomorphic to a complemented subspace of $(\sum
H^1_n)_{l^1}$.
\item If $$\sup_{Q\in\cE}\sup_{P\sbe Q, P\in \cE}\mu(P)/\mu(Q)<\infty,$$
then for $1\le j\le N$ $H^1([\cF^j_n]^\infty_{n=1})$ is isomorphic  to a
complemented subspace of $l^1$. \end{enumerate}

\section{A smooth unconditional basic sequence in $L^2(X,\mu)$}

Let $Q$ be  in $G_n(X|\cE)$ and let $P_0,P_1,\dots,P_N$ be an enumeration of
$G_1(Q|\cE)$. By the above preliminary remarks, $N\le N_0$, where $N_0$ is
independent of $Q$ and
$$\inf_i\frac{\mu(P_i)}{\mu(Q)}>\frac 1C$$
$$\frac
1C\le\inf_i\frac{\mu(P_i)}{\mu(P_0)}\le\sup_i\frac{\mu(P_i)}{\mu(P_0)}\le C$$
where $C$, depending on $\d$, is independent of $Q$.

Let, for $1<i\le N$, the function $h_{Q,i}$ satisfy the following conditions
\begin{enumerate}
\item $\supp h_{Q,i}\sb Q$
\item $h_{Q,i}$ is constant when restricted to one the sets $P_j,0\le j\le N$.
\item There exists $C>0$ (not depending on $\d$ or $N$) so that for $\a_i\in
\tR$
$$\frac 1C\left(\sum^N\a^2_i\right)\le
\biggl\Vert\sum^N_{i=1}h_{Q,i}\a_i\biggr\Vert_{L_2(X,\mu)}\le
C\left(\sum\a^2_i\right)^{1/2}.$$
\end{enumerate}
Using ideas related to the local Pe\l czy\'nski decomposition, such a stystem
was constructed by B. Maurey [Ma1].

As martingale differences an orthogonal in $L^2(X,\mu)$ we get for $f\in
L^2(X,\mu)$ a uniquely determined sequence of coefficients $\a_{Q,i}$,
$Q\in\cE$ so that
$$f=\sum_{Q\in\cE}\sum_{i\in I_Q}h_{Q,i}\a_{Q,i}$$
and
\begin{eqnarray*}||f||_2&=&\left(\sum_{Q\in\cE}\biggl\Vert\sum_{i\in I_Q}
h_{Q,i},\a_{Q,i}\biggr\Vert^2_2\right)^{1/2}\\
&\sim&\left(\sum_{Q\in\cE}|Q|\sum_{i\in
I_Q}\a^2_{Q,i}\right)^{1/2}.\end{eqnarray*}
In the other words $h_{Q,i}$  $Q\in\cE$, $i\in I_Q$ forms an unconditional basis
in $L^2(X,\mu)$. Using smoth partition of unity we will modify $h_{Q_i}$ to
become a smoth unconditional basis for
$L^2(X,\mu)$.

For $Q\in G_n(X|\cE)$ we have
$$C_2\d^n\le\diam Q<C_1\d^n.$$
For $\tau<1/C_1500$ we consider a partition of
unity $\psi_k^{(n)}, k=1,\dots,N_n$,
so that:
\begin{eqnarray*} \diam(\supp\psi_k^n)&\le& \tau\d^n\\
\Lip_\b(\psi^n_k)&\le&(\tau\d^n)^{-\b}\\
\sum^{N_n}_{k=1} \psi^{(n)}&=&1.\end{eqnarray*}

See [M.-S.2] for a construction of such a partition of unity. We use it here to
define the kernel
$$K_n(x,y):=\sum^{N_n}_{k=1}\psi^{(n)}_k(x)\psi^{(n)}_k(y)\frac
1{||\psi^{(n)}_k||_1}$$
and define
\begin{eqnarray*}
\tilde\vp_{Q,i}(x)&:=&\int_XK_{n+1}(x,y)h_{Q,i}(y)d\mu(y)\\
\vp_{Q,i}(x)&:=&\frac{\tilde\vp_{Q,i}(x)}{||\vp_{Q,i}||_2}.\end{eqnarray*}
                                    *
By construction we obtain at once the following properties of $\vp_{Q,i}$:
\begin{eqnarray*} \supp\vp_{Q,i}&\sb&\{z\in X,\dist(\supp
h_{Q,i},z)\le\tau\d^n\}\\
\Lip_\b(\vp_{Q,i})&\le&\left(\frac{\mu(Q)}{\tau\d}\right)^\b\left(\frac{\mu(Q)}
{\d}\right)^{-1/2}\\
\int_X\vp_{Q,i}d\mu&=&0.\end{eqnarray*}
And for $Q\in\cE$ fixed we obtain
$$\frac
1C\left(\sum_i\a^2_i\right)^{1/2}\le\biggl\Vert\sum_i\a_i\vp_{Q,i}
\biggr\Vert_{L^2(X,\mu)}\le
\left(\sum_i\a^2_i\right)^{1/2}C$$
where $C$ is independent of $Q$ or $\tau$ and  depends only on the geometry of
$(X,d,\mu)$.

Moreover we have the following theorem.

\begin{theor} Let $\cE_1=\bigcup^\infty_{n=1}G_{2n}(X|\cE)$ then for
$$f=\sum_{Q\in\cE_1}\sum_{i\in IQ}\a_{Q,i}\vp_{Q,i}$$
we have
$$\left(\sum_Q\sum\a^2_{Q,i}\right)^{1/2}\le||f||_2\le
\left(\sum_Q\sum\a^2_{Q,i}\right)^{1/2}.$$
       \end{theor}

\proof Suppose $\supp \vp_{Q,i}\cap\supp\vp_{Q,i}\ne 0$ then w.l.o.g. assume that
$$\diam(\supp\vp_{Q,i})\le \diam(\supp \vp_{P,j}).$$
Let $z$ be a fixed point in $\supp \vp_{Q,i}$, then
\begin{eqnarray*}
\int_X\vp_{Q,i}\vp_{P,j}d\mu&=&\int_X\vp_{Q,i}(\vp_{p,j}-\vp_{p,j}(z))d\mu\\
&\le&||\vp_{Q,i}||_1\sup_{x\in Q}|\vp_{p,j}(x)-\vp_{p,j}(z)
|\le |\d^{1/2}\mu(Q)^{1/2}(\Lip_\b\vp_{p,j})\diam
(\supp \vp_Q)^\b|\\
&\le&\frac{\mu(Q)^{1/2+\b}} {\mu(P)^{1/2+\b}}\frac
1{(\tau\d)^\b}.\end{eqnarray*}

Then given $Q$, consider $P\in G_n(Q|\cE_1)$ then $\frac{\mu(P)}{\mu(Q)}\le
C\d^{2n}$ and $G_n(Q|\cE)$ contains at most $C\d^{-2n}$ elements.

From these observations we see (using e.g. the argument in [U, Lemma 3.3]) that
there exists $C$ (not depending on $\d$ or $\tau$) so that
$$C||f||^2_2+\frac{\d^{2\b}}{(\tau\d)^\b}\sum\sum\a^2_{Q,i}\ge\frac
1C\sum\sum\a^2_{Q,i}$$
and
$$||f||^2_2<2C\sum\sum\a^2_{Q,i}.$$

Now choosing $\d$ so small that $(\frac \d \tau)^\b<\frac 1{C^2}$ we
obtain the result.

\section{A smooth biorthogonal sequence in $L^2(X,\mu)$}

Let $G_n:=G_n(X|\cE)$. Fix $K\gg 1$. Using Lemma 9 from [M\"u3] we split
$G_n$ into
$\cP_{n,1},\dots,\cP_{n,l}$ so that for
$P,Q\in \cP_{n,j}$ we have $\dist(P,Q)\ge K\mu\{\mu(P),\mu(Q)\}$ and $l$ depends
only on $K$ and the geometry of $(X,d,\mu)$.

Now fix $m\in\tN\sm \{1\}$ $0<s\le m$ and $j\le l$.
Then let $\cF:=\bigcup^\infty_{k=0}\cP_{mk +s,j}\cup\{X\}$. Next fix $i_0\le N$
and for $Q\in\cF$ let
\begin{eqnarray*} \vp_Q&:=&\vp_{Q,i_0}\\
\vp_X&:=&\1_X.\end{eqnarray*}

Observe now that for $P,Q\in P_{m\circ k+s,j}$
$$\supp\vp_Q\cap\supp\vp_P\ne\es$$
and for each $P\in \cP_{mk+s,j}$ and $r\in\tN$ there exists at most one $Q\in
\cP_{m(k-r)+s,j}$ so that
$$\supp\vp_Q\cap \supp\vp_P=\es.$$

 Moreover by Theorem 3 the Gram matrix
$$G:=\left(\int \vp_Q\vp_Pd\mu\right)_{Q,P\in\cF}$$
is invertible (and positive definite).

The Gram-matrix is used to construct a biorthogonal system from the $\vp_Q-s$.

\begin{theor}
\begin{itemize}
\item[a)] The coefficients $(a_{P,Q})_{P,Q\in\cF}$ of the matrix $G^{-1/2}$
satisfy the estimates
$$|a_{P,Q}|\le
C\min\left\{\frac{\mu(P)}{\mu(Q)}\frac{\mu(Q)},{\mu(P)}\right\}^{1/2-
\a}\left(1+\frac{\dist (P,Q)}{3\max\{\mu(P),\mu(Q)\}}\right)^{-1-\a}$$
where $0<\a<\b/2$.
\item[b)] The functions $$f_Q:=\sum_{P\in\cF}a_{P,Q}\vp_P,\quad\quad\quad
Q\in\cF$$
form an orthonormal system in $L^2(X,\mu)$, the closed span of which coincides
with the closed span of $\{\vp_Q:Q\in\cF\}$.\end{itemize} \end{theor}

\proof Part b) is a well known algebraic identity, so we shall concentrate on
the Proof of part a):

Recall first that for $Q\cap P\ne 0$ we have the estimate
$$\int\vp_Q\vp_Pd\mu\le\left(\frac
1{\tau\d}\right)^\b\min\left\{\frac{\mu(Q)}{\mu(P)},\frac{\mu(P)}{\mu(Q)}
\right\}^{1/2+\b}.$$

Combining this with $\supp\vp_Q\sbe \{z\in X:d(z,Q)\le\a\}$ we obtain in
particular
$$\int\vp_Q\vp_Pd\mu\le \min\left\{\frac{\mu(P)}{\mu(Q)},\frac
{\mu(Q)}{\mu(P)}\right\}^{1/2+\b}\left(1+\frac{\dist(P,Q)}{\max\{\mu(P),\mu(Q)\}}
\right)^{-1-\b}.\leqno(3.1)$$

Moreover if $Q,P\in\cF$, $Q\ne P$ and $\int \vp_Q\vp_P\ne 0$ then necessarily
$$\min  \left\{\frac{\mu(P)}{\mu(Q)},\frac
{\mu(Q)}{\mu(P)}\right\}\le \d^m.$$
(To obtain this conclusion we introduced the splitting of $G_n$ into $\cP_{n,j}$.)

Using this information, the proof of Lemma 3.3 in [U] gives that for
$a:=(a_P)_{P\in\cF}$, $a_P\in\tR$ the following norm estimate for the matrix
$\Id-G$ holds:

$$||(\Id-G)a||_{l^2}\le \left(\frac{\d^m}{(\d\tau)^\b}\right)C_2||a||_{l^2}.$$
($C_2$ is a universal constant.)
Now put $R=\Id-G$. Observe that $G^{-1/2}$ can be developed in a power series of
$R$, indeed:
$$G^{-1/2}=\sum^\infty_{k=0}C_kR^k$$
where $C_k=o(k^{-1/2})$.

Clearly the coefficients $R(P,Q)$, $P,Q\in\cF$ of $R$ satisfy the estimates
(3.1).

Now by a result of Frazier-Jawerth, estimates of this form are stable under the
formation of products. More precisely by [F-J, Theorem 9.1] for $0<\g<\b$ there
exists $C_1>1$
so that for each $k\in\tN$ the coefficients $R^{(k)}(P,Q)$ of
$R^k$ satisfy
$$R^{(k)}(P,Q)\le\left[\left(\frac
1{\tau\b}\right)^\b C_1\right]^k\min\left\{\frac{\mu(P)}{\mu(Q)},\frac
{\mu(Q)}{\mu(P)}\right\} ^{1/2+\g}$$
$$\left(1+\frac{\dist(P,Q)}{\max\{\mu(P),\mu(Q)\}}\right)^{-1-\g}.$$
On the other hand we trivially have
$$R^{(k)}(P,Q)\le||R^k||_{l^2}\le\left[\frac{\d^mC_2}{(\tau\d)^\b}\right]^k.$$

Fix now $P,Q\in\cF$ and let
$$\s(P,Q):=\left\{\min\frac{\mu(P)}{\mu(Q)},\frac
{\mu(Q)}{\mu(P)}\right\}^{1/2+\g} \left(1+\frac{\dist(P,Q)}{\max\{\mu(P),
\mu(Q)\}}\right)^{-1-\g}.$$

Next consider the number $k_0=k_0(P,Q)$ which is defined by
$$k_0:=\left[\frac\g 2\frac{\log\s(P,Q)}{\log(C_1/(\d \tau)^\b)}\right].$$
We assume that $k_0$ is integer.
At this point we make a suitable choice for $m$. Namely we choose $m$ so that
$$\frac{\log(C_2\d^m/(\d  \tau)^\b)}{\log ((\tau\d)^\b C_1)}\ge 1+\frac \g 2.$$
(Observe that $m$ is of course {\bf not} depending on $P,Q$.)

We then have the numerical estimates:
$$\sum^{k_0}_{k=0}\left(\frac{C_1}{(\d\tau)^\b}\right)^k\le\s(P,Q)^{-\g/2}C$$
$$\sum^\infty_{k=k_0+1}\left(\frac{\d^m}{(\tau\d)^\b}C_2 \right)
^k\le\s(P,Q)^{1+\g/2}C.$$

Hence the following estimates hold for $P,Q\in\cF$
$$\sum^{k_0}_{k=0}R^{(k)}(P,Q)\le \s(P,Q)^{1+\g-\g/2}C$$
$$\sum^\infty_{k=k_0+1}R^{(k)}(P,Q)\le\s(P,Q)^{1+\g/2}C.$$

Summing up we have for the coefficients of $G^{-1/2}$ the following estimate
$$G^{-1/2}(P,Q)\le C\min \left\{\frac{\mu(P)}{\mu(Q)},\frac
{\mu(Q)}{\mu(P)}\right\} ^{1/2+\g/2}
\left(1+\frac{\dist(P,Q)}{\max\{\mu(P),
\mu(Q)\}}\right)^{-1-\g/2}C.$$

\rem The above proof merges arguments of Franzier \& Jawerth [F,J] and Uchiyama
[U] with those of Jaffard \& Meyer [J-M] to conclude that -- in the language if
Franzier and Jawerth -- $G^{-1/2}$ is an almost diagonal matrix. As a
consequence $f_Q$ is centered around $Q$. More precisely we have the following
pointwise estimate.
\endproof

\begin{lemma} There exists $C=C(\d)\sim\log\d$ so that
\begin{enumerate}
\item for $x\in X$
$$|f_Q(x)|\le C\left(1+\frac{\dist(x,Q)}{\mu(Q)}\right)^{-1-\a/2}\frac
1{\mu(Q)^{1/2}}$$
\item for $x,y\in X$ with $d(x,y)\le\mu(Q)$
$$|f_Q(x)-f_Q(y)|\le C\left(1+\frac{\dist(x,Q)}{\mu(Q)}\right)^{-1-\a/2-\b}
\frac{d(x,y)^\b}{\mu(Q)^{1/2+\b}}.$$
\end{enumerate}\end{lemma}

\proof Given the estimates of $G^{-1/2}$ the proof is quite standard, and only
the argument for part 1) will be outlined. Fix $x\in X$ then clearly for some
$C>1$
\begin{eqnarray*}|f_Q(x)|&\le&C\sum_{\{P:\dist(x,P)\le C\mu(P)\}}\min
\left\{\frac{\mu(P)}{\mu(Q)},\frac
{\mu(Q)}{\mu(P)}\right\} ^{1/2+\a}\\
&\times&\frac 1{\mu(P)^{1/2}}\left(1+\frac{\dist(P,Q)}{\max\{\mu(P)\mu(Q)}\right)
^{-1-\a}.\end{eqnarray*}
\endproof

Let now $K\in\tN$ be such that
$$K\mu(Q)\le\dist(x,Q)\le 2K\mu(Q),$$
and split the above sum into three, by dividing the index set:
$$\{P:\dist(x,P)<C\mu(P)\}=\cA\cup \cB\cup \cD$$ where
$$\cA:=\{P\dist(x,P)\le C\mu(P)\mbox{ and } \mu(P)>K\mu(Q)\}$$
$$\cB:=\{P:\dist(x,P)\le C\mu(P)\mbox{ and } \mu(Q)\le \mu(P)\le K\mu(Q)\}$$
$$\cD:=\{P:=dist(x,P)\le C\mu(P)\mbox{ and } \mu(P)\le mu(Q)\}.$$

{\bf Case 1}
$$\sum_{\cA}\left\{\frac{\mu(Q)}{\mu(P)}\right\}^{1/2+\a}
\left(1+\frac{\dist(P,Q)}
{\mu(P)}\right)
^{-1-\a}\frac 1{\mu(P)^{1/2}}$$
$$\le\sum_{\mu(P)>K\mu(Q)}\left\{\frac{\mu(Q)}{\mu(P)}\right\}^{1/2+\a} \frac
1{\mu(P)^{1/2}}$$
$$\le\frac{\mu(Q)^{\a+1/2}}{(K\mu(Q))^{1+\a}}=\frac 1{\mu(Q)^{1/2}K^{1+\a}}$$
$$\le \frac 1{\mu(Q)^{1/2}}\left(1+\frac{\dist(x,Q)}{\mu(Q)}\right)^{-1-\a}$$

{\bf Case 2} Let $\mu(P)=\mu(Q)\d^{-k}$ and $\d^{-k_1}=K$.
Then

$$\sum_{\cB}\left\{\frac{\mu(Q)}{\mu(P)}\right\}^{1/2+\a}
\left(1+\frac{\dist(P,Q)}
{\mu(P)}\right)
^{-1-\a}\frac 1{\mu(P)^{1/2}}$$

\begin{eqnarray*}
&\le& \sum^{k_1}_{k=0} \d^{k(1/2+\a)}\d^{k/2}\left(\frac{\mu(Q)\d^k}{k}\right)
^{-1-\a}\frac 1{\mu(Q)^{1/2}}\\
&\le&\mu(Q)^{1+\a}\mu(Q)^{-1/2}\frac{k_1}{K^{1+\a}}\\
&\le&\mu(Q)^{-1/2}|\log \d|\frac{\log K}{K^{1+\a}}\\
&\le&c_\a(\log\d)\mu(Q)^{-1/2}\left(1+\frac{\dist(x,Q)}{\mu(Q)}\right)
^{-1-\a/2}\end{eqnarray*}

{\bf Case 3}
$$\sum_{\cD}\left\{\frac{\mu(Q)}{\mu(P)}\right\}^{1/2+\a}
\left(1+\frac{d(P,Q)}{\mu(Q)}\right)^{-1-\a}\frac
1{\mu(P)^{1/2}}\le\left\{\begin{array}{ll} \frac 1{\mu(Q)^{1/2}}
&\mbox{if }\dist (x,Q)\le C\mu(Q)\\
0&\mbox{otherwise.}\end{array}\right.$$

\section{Bounded Projections in $H^1(X,d,\mu)$}

First we determine the norm of $f_Q,Q\in\cF$ in $H^1(X,d,\mu)$. Given the decay
of $f_Q$ and the fact that $\int_Xf_Qd\mu=0$ it is natural to use molecules as
in [Woj].

\begin{theor}  There exists $C=C(\d,\a)$ and $\ve>0$ so that for each $Q\in\cF$
we have
$$\left(\int f^2_Q\frac{d\mu}{\mu(Q)}\right)\left(\int
f^2_Q(x)d(x,Q)^{1-\ve}\frac{d\mu}{\mu(Q)}\right)^{1/\ve}\le C(\d,\a).$$\end{theor}

\proof First we have clearly
$$\int f^2_Q\frac{d\mu}{\mu(Q)}=\frac 1{\mu(Q)}.$$
Let $Q_n:=\{x\in X,\mu(Q)(2^n-1)\le d(x,Q)\le(2^{n+1}-1)\mu(Q)\}$,
then
$$\int
f^2_Q(x)d(x,Q)^{1+\ve}d\mu=\sum^\infty_{n=0}\int_{Q_n}f^2_Q(x)d(x,Q)
^{1+\ve}d\mu.$$
Let us first consider the case $n\ge 1$:

\begin{eqnarray*} \int_{Q_n}f^2_Q(x)d(x,Q)^{1+\ve} d\mu&\le&
C(\d)\int_{Q_n}\left(1+\frac{d(x,Q)}{\mu(Q)}\right)^{-2-\a}(d(x,Q)^{1+\ve}
\frac{d\mu}{\mu(Q)^{1/2}}\\
&\le&C(\d) 2^{n(-2-\a)}\mu(Q)^{-1}(2^{n+1}\mu(Q))^{1+\ve}2^{n+1}\mu(Q)\\
&\le&C(\d)\mu(Q)^{1+\ve}4\cdot 2^{n(-\a+\ve)}.\end{eqnarray*}

And for $n=0$  we have
$$\int_{Q_0}f^2_Q(x)d(x,Q)^{1+\ve}d\mu\le \mu\frac
1{\mu(Q)}\int_{Q_0}d(x,Q)^{1+\ve}d\mu\le\mu(Q)^{1+\ve}.$$
Summing up we obtain
$$\int_Xf_Q^2(x)d(x,Q)^{1+\ve}d\mu\le C(\d)\frac 1{1-2^{-\a+\ve}}\mu(Q)$$
and finally
$$\left(\int f^2_Q\frac{d\mu}{\mu(Q)} \right)\left(\int
f^2_Q(x)d(x,Q)^{1+\ve}\frac{d\mu}{\mu(Q)}\right)^{1/ \ve}\le\left(C(\d)
\frac 1{1-2^{-\a+\ve}}\right)^{1/ \ve}.$$
Choosing $\ve=\a/2$ gives the required estimate.

We shall show next that $\ol{\span}\{f_Q:Q\in\cF\}$ is a complemented subspace
of $H^1(X,d,\mu)$ and that $\{f_Q:Q\in\cF\}$ is equivalent to a martingale
difference sequence in $H^1([\cF_n])$ (where $\cF_n$ was defined  in Section 1).

The operator $Pf=\sum_{Q\in\cF}(f|f_Q)f_Q$ is clearly a projection , i.e.,
satisfies\\ $P^2=P$.

Theorem together with the smothness and localization properties of $f_Q$
will be used to show that $P$ defines a bounded projection on $H^1(X,d,\mu)$.

\begin{theor} There exists $C>0$ so that for $f\in H^1(X,d,\mu)$
$$\biggl\Vert
\sum_{Q\in\cF}(f|f_Q)f_Q\biggr\Vert_{H^1(X,d,\mu)}\le||f||_{H^1(X,d,\mu)}.$$
\end{theor}

\rem The following proof not new! It is a simple modification of the proof in
[Woj Theorem], and is included here just for sake of completeness.
\bs

\proof It es enough to consider atoms in $(X,d,\mu)$: Let $a:X\to\tR$ be
supported on a ball $B$ so that $\int ad\mu=0$, $||a||_\infty\le\mu(B)^{-1}$ and
$\mu(B)\le C\diam B$.
Then decompose $\cF=E\cup F\cup G$ where

\begin{eqnarray*}
E&=&\{Q\in\cF:\mu(Q)\ge\mu(B)\}\\
F&=&\{Q\in\cF:\mu(Q)\ge\mu(B)\mbox{ and } \dist(P,Q)\le L\mu(Q)\}\\
G&=&\{Q\in\cF:\mu(Q)\ge\mu(B) \mbox{ and } \dist(P,Q)\ge
L\mu(Q)\}.\end{eqnarray*}\endproof

{\bf Case 1} By the  triangle inequality we have: using Theorem 6:

\begin{eqnarray*} \biggl\Vert\sum_{Q\in E}(a|f_Q)f_Q\biggr\Vert_{H^1}&\le&
\sum_{Q\in E} |(a|f_Q)|\cdot ||f_Q||_{H^1}\\
&\le&\sum_{Q\in E}\mu(Q)^{1/2}\frac{\diam
B^\b}{\mu(Q)^{1/2+\b}}\left\{1+\frac{d(B,Q)}{\mu(Q)}\right\}^{-1-\a/2-\b}\\
&\le&\sum_{\mu(Q)>\mu(B)}\left\{\sum^{\mu(Q)^{-1}}_{k=1}k^{-1-\a/2-\b}
\right\}\diam
B^\b\mu(Q)^{-\b}\\
&\le&\frac 1{\a/2+\b}\left\{\sum_{\mu(Q)>\mu(B)}\mu (Q)^{-\b}\right\}\diam (B)\\
&\le&\frac 1{\a/2+\b}\mu(B)^{-\b}\diam (B)\le\mbox{ \rm const.}\end{eqnarray*}

{\bf Case 2} Again by triangle inequality and Theorem 6:

\begin{eqnarray*}
\biggl\Vert \sum_{Q\in G} (a|f_Q)f_Q\biggr\Vert_{H^1}&\le& \sum_{Q\in
G}\int_B|f_Q|d\mu\mu(B)^{-1}\mu(Q)^{1/2}\\
&\le&C\sum_{Q\in
G}\mu(Q)^{-1/2}\left(1+\frac{\dist(B,Q)}{\mu(Q)}\right)^{-1-\a/2}\mu(Q)^{1/2}\\
&\le&  C\sum_{\mu(Q)\le\mu(B)}\left(1+\frac{\mu(B)}{\mu(Q)}\right)^{-\a/2}\\
&\le&C\sum_{\mu(Q)\le\mu(B)}\left(\frac{\mu(Q)}{\mu(B)}\right)^{\a/2}\le
C\mbox{ \rm const.}\end{eqnarray*}

{\bf Case 3} Here we show that
$\sum_{Q\in F}(a|f_Q)f_Q$ is a molecule.

Consider first
$$\int_B\biggl\Vert\sum_{Q\in F} (a|f_Q)f_Q\biggr\Vert^2d(x,x_B)^{1+\ve}\le
C\mu(B)^{1+\ve}||a||^2_2\le c\mu(B)^\ve.$$
Then we consider
\begin{eqnarray*}
\lefteqn{\int_{X\sm B}\biggl\Vert\sum_{Q\in F}(a|f_Q)f_Q\biggr\Vert^2d(x,x_B)
^{1+\ve}d\mu(x)}\\
&&\le  
 C||a||^2_2\sum_{Q\in F}\mu(Q)^{-1}\int_{X\sm
B}\left(\frac{1+d(x,Q)}{\mu(Q)}\right)^{-2-\a}d(x,x_B)^{1+\ve}d\mu\\
&&\le C\mu(B)^{-1}\sum_{\mu(Q)\le\mu(B)\mu(B)^{-1}}\frac{\mu(B)}{\mu(Q)}
\int_{X\sm
B}\left(\frac{d(x,x_B)}{\mu(Q)}\right)^{-2-\1}d(x,x_B)^{1+\ve}d\mu\\
&&\le C\sum_{\mu(Q)\le\mu(B)}\mu(Q)^\a\int_{X\sm B}d(x,x_B)^{-1+\ve-\a}d\mu\\
&&\le C\left\{\sum_{\mu(Q)\le\mu B)} \mu(Q)^\a \mu(B)^\a\right\}\mu(B)^\ve.
\end{eqnarray*}

Summing up we have for $\ve<\a$:
$$\left(\int_X\biggl\Vert\sum_{Q\in
F}(a|f_Q)f_Q\biggr\Vert^2d(x,x_B)^{1+\ve}\right)^{1/ \ve}\le C\mu(B)$$
and
$$\int_X \biggl\Vert\sum_{Q\in F}(a|f_Q)f_Q\biggr\Vert^2d\mu\le||a||^2_2\le
C\mu(B)^{-1}.$$
Multiplying the above estimates one sees that
$\sum_{Q\in F}(a|f_Q)f_Q$ is indeed a molecule.
\endproof

In Section 1, using successive, generations of $\ve$, an increasing sequence of
$\a$-algebra, $(\cF_n)^\infty_{n=1}$ has been defined.

In Section 2, we defined on unconditional basis $\{h_{Q,i},Q\in \ve,i\in I_Q\}$
for $L^2(X,\mu)$. As recorded in [Ma2] this system forms an unconditional basis
in the martingale $H^1([\cF_n])$ space.

We fix now $i_0\le N$ as in Section 3 and let
$$h_Q=h_{Q,i_0},\quad\quad\quad Q\in\cF.$$
The family $\{h_Q:Q\in\cF\}$ forms a three valued martingale difference sequence
with respect to the filtration $[\cF_n]^\infty_{n=1}$ satisfying the following
condition:
$$\supp h_Q\cap \supp h_P\ne\es$$ implies
$$\supp h_Q\sbe\supp h_P\mbox{ or } \supp h_P\sbe \supp h_Q.$$

We will show next, that $\{f_Q,Q\in\cF\}$ in $H^1(X,d,\mu)$ is equivalent to
$\{h_Q:Q\in\cF\}$ in $H^1([\cF_n]\}$.

Let $Y$ be the closed linear span of $\{f_Q:Q\in\cF\}$ equipped with the norm
inherited by $H^1(X,d,\mu)$, then we have:

\begin{theor} \begin{eqnarray*} T:Y&\to&H^1([\cF_n])\\
f_Q&\to&h_Q\end{eqnarray*}
extendes to a bounded operator.
\end{theor}

\proof Let $f\in Y$ implies clearly $f\in H^1(X,d,\mu)$. Hence there exist
atoms $a_i$, and $\lambda_i\in\tR$ so that
$$f=\sum\lambda_ia_i\mbox{ and } \sum|\lambda_i|\le C||f||_{H^1}.$$
Moreover
$$f=Pf=\sum\lambda_iPa_i$$
and
$$||Pa_i||_{H^1(x,d,\mu)}\le C||a_i||_{H^1(X,d,\mu)}.$$
So it remains to show that there exists $C>0$ so that for any atom $a$ on
$(X,d,\mu)$ we have $||TPa||_{H^1([\cF_n])}\le C$.

To estimate
$$TP_a=\sum_{Q\in\cF}(a|f_Q)h_Q$$
in $H^1([\cF_n])$ we observe that
$$||h_Q||_{H^1}([\cF_n])\le C\mu(Q)^{1/2},$$
split $\cF$ into $E\cup F\cup G$ as in the proof of Theorem (7) and argue
exactly as P. Wojtaszczyk in [Woj2, Theorem 5].\endproof

Let $Z$ be closure of the linear span of $\{h_Q:Q\in\cF\}$ in $H^1([\cF_n])$,
equipped with the norm inherited by $H^1([\cF_n])$.
By [Ma1, Theorem 2] $\{h_{Q,i},Q\in\ve,i\in
I_Q\}$ is an {\bf unconditional} basis in $H^1([\cF_n])$ the natural restriction
operator
\begin{eqnarray*}Q:H^1([\cF_n])&\to& H^1([\cF_n])\\
\sum_{Q\in\ve}\sum_{i\in I_Q}\a_{Q,i}h_{Q,i}&\to&
\sum_{Q\in\cF}a_{Q,i_0}h_{Q,i_0}\end{eqnarray*}
is a bounded projection.

Moreover given any atom $a$ in the martingale $H^1([\cF_n])$ space then $Qa$
is again an atom in $H^1([\cF_n])$.

(In Section 1 we remarked that the filtration $[\cF_n]_{n=1}^\infty$ is regular
(see [G, p 96]) and therefore an atom in $H^1([\cF_n])$ is simply a function
$a:X\to \tR$ for which is supported in an atom $Q$ of $\cF_n$ so that
$||a||_\infty\le \mu(Q)^{-1}C$ and $\int ad\mu=0$.)

Now we have the following

\begin{theor} \begin{eqnarray*} S:Z&\to&H^1(X,d,\mu)\\
h_Q&\to&f_Q\end{eqnarray*}
defines a bounded operator.
\end{theor}

\proof Let $f\in Z$. Then there exists a sequence of atoms $a_i$ for
$H^1([\cF_n])$ and $\lambda_i\in\tR$ so that
$$f=\sum\lambda_ia_i$$
and
$$\sum^\infty_{i=1}|\lambda_i|\le C||f||_{H^1}([\cF_n]).$$
As $$f=Qf=\sum^\infty_{i=1} \lambda_iQa_i,$$

we have: that for  any $f\in Z$ there exists a sequence of atoms $q_i$: in
$H^1([\cF_n])$, $\lambda_i\in\tR$ and $q_i\in Z$ (sic!) satisfying
$$f=\sum^\infty_{i=1}\lambda_iq_i\mbox{ and }\sum|\lambda_i|\le
C||f||_{H^1([\cF_n])}.$$
It is therefore enough to consider atoms $q$ of the form
$$q=\sum_{Q\in\cF}\a_Qh_Q$$
and to show that
$$||Sq||_{H^1(X,d,\mu)}=\biggl\Vert
\sum_{Q\in\cF}\a_Qf_Q\biggr\Vert_{H^1(X,d,\mu)}$$
is bounded by an absolute constant independent of $q$.

As moreover $\{h_Q:Q\in\cF\}$ is biorthogonal it remains to show that  there
exists $C>0$ so that for any atom $q\in Z$
$$\biggl\Vert \sum_{Q\in\cF}(q|h_Q)f_Q\biggr\Vert_{H^1(X,d,\mu)}\le C.$$
To do so we just follow the argument in [Woj Theorem 5] again.
            \endproof

\section {D\'enoument}
In this paragraph we will give a solution to the
classification problem of atomic $H^1(X,d,\mu)$ spaces:

In addition to the material developed in Sections 1 --- 4 we will use the
following ingredients
\begin{itemize}
\item The isomorphic classification of martingal $H^1$-spaces generated by an
increasing sequence of purely atomic $\s$-algebras.
\item The isomorphic classification three-valued martingale difference sequences
in martingale $H^1$ spaces.
\item $H^1(X,d,\mu)$ is isomorphic to a complemented subspace of martingale
$H^1$ space.
      \end{itemize}

\begin{theor} If $H^1(X,d,\mu)$ is infinite dimensional, it is isomorphic to one
of the following spaces: $H^1(\d), (\sum H^1_n)_{l^1},l^1$.\end{theor}

\proof \begin{enumerate}
\item {\bf The Case $H^1(\d)$}

Let $E=\{t\in X:t$ lies in infinitely many elements of $\ve\}$. Suppose
$\mu(E)>0$. Then there exists a subcollection $\cF\sb\ve$ as constructed in
Section 3 so that
$$F:=\{t\in X:t\mbox{ lies in infinitely many elements of } \cF\}$$
satisfies $\mu(F)>0$.

By [Mu2], $\ol{\span}\{h_Q:Q\in\cF\}$ equipped with the norm of
$H^1([\cF_n])$ is
then isomorphic to $H^1(\d)$. Hence by Section 4
$$H^1(\d)\stackrel{C}{\hookrightarrow} H^1(X,d,\mu).$$
On the other hand by the results in [M\"u3] and [Ma3]
$$H^1(X,d,\mu)\stackrel{C}{\hookrightarrow} H^1(\d).$$
So the Pe\l czy\'nski decomposition method gives that $H^1(\d)$ is isomorphic to
$H^1(X,d,\mu)$.

\item  {\bf The Case $(\sum H^1_n)_{l^1}$}

Suppose that $\mu(E)=0$ and
$\sup_{Q\in \ve} \sum_{P\sb Q,P\in \ve}\mu(P)/ \mu (Q)=\infty$.
Then there exists a subcollection $\cF\sb\ve$ constructed as in Section 3 so
that $\mu(F)=0$ and
$$\sup_{Q\in \cF}\sum_{P\sb Q,P\in\cF}\mu(P)/ \mu (Q)=\infty.$$

By the result of [Mu2] $\ol{\span}\{h_Q:Q\in\cF\}$ is then isomorphic to
$(\sum|H^1_n)_{l^1}$. Hence by Section 4
$$(\sum H^1_n)_{l^1}\stackrel{C}{\hookrightarrow} H^1(X,d,\mu).$$
On the other hand by [Mu2] $\mu(E)=0$ implies
$$H^1(X,d,\mu)\stackrel{C}{\hookrightarrow} \left(\sum H^1_n\right)_{l^1}.$$
So the Pe\l czy\'nski decomposition method gives   that $(\sum H^1_n)_{l^1}$ is
isomorphic to $H^1(X,d,\mu)$.

\item {\bf The Case $l^1$}

Suppose
$$\sup _{Q\in\ve}\sum_{P\sb Q,P\in\ve}\mu(P)/ \mu(Q)<\infty$$
then by [M\"u1] and [M\"u3]
$$H^1(X,d,\,u)\stackrel{C}{\hookrightarrow} l^1.$$
By a theorem of Pe\l czy\'nski a complemented subspace of $l^1$ is either
finite dimensional or isomorhic to $l^1$.

\end{enumerate}
\endproof

\vskip 1.0cm
Institut f\"ur Mathematik\\
Johannes Kepler Universit\"at\\
A-4040 Linz\\
Austria

\end{document}